\documentclass{elsart}
\usepackage{amssymb}
\usepackage{array}
\usepackage{fancyhdr}
\usepackage{longtable}
\usepackage{tabularx}
\usepackage{multirow}
\usepackage{rotating}

\usepackage{amsmath}
\usepackage{theorem}
\usepackage[mathscr]{eucal}
\usepackage{enumerate}
\usepackage[dvips]{epsfig}
\usepackage{float}  
\journal{Applied Numerical Mathematics}
\volume{59}
\pubyear{2009}
\newcounter{mylastpage}
\issue{3--4}
\usepackage{hyperref}
\makeatletter
\def\ps@copyright{%
 \let\@oddhead\@empty
 \let\@evenhead\@empty
 \def\@oddfoot{\small
 \slshape\hskip-5em
   Published in \@journal\ \@volume\ (\the\@pubyear)\ no.\ \@issue, pp.\ \ESpagenumber{firstpage}--\ESpagenumber{mylastpage},\newline
\href{http://dx.doi.org/10.1016/j.apnum.2008.03.011}{doi: 10.1016/j.apnum.2008.03.011}\hskip-70pt}%
 \let\@evenfoot\@oddfoot}
\makeatother
\newtheorem{Def}{Definition}[section]
\newtheorem{Sat}[Def]{Proposition}

\newtheorem{Lem}[Def]{Lemma}

\newcommand{\arsinh}{\operatorname{arsinh}}

\newcommand{\E}{\operatorname{E}}
\newcommand{\Prob}{\operatorname{P}}

\begin{document}
\begin{frontmatter}
\title{Diagonally Drift--Implicit Runge--Kutta Methods of Weak Order One and Two for It\^{o} SDEs
and Stability Analysis}
\author{Kristian Debrabant} and
\ead{debrabant@mathematik.tu-darmstadt.de}
\author{Andreas R\"{o}{\ss}ler}
\ead{roessler@mathematik.tu-darmstadt.de}
\address{Technische Universit\"{a}t Darmstadt, Fachbereich Mathematik, Schlossgartenstr.7,
D-64289 Darmstadt, Germany}
{\em In honor of Professor Karl Strehmel}
\begin{abstract}
The class of stochastic Runge--Kutta methods for stochastic
differential equations due to R\"{o}{\ss}ler is considered. Coefficient
families of diagonally drift--implicit stochastic Runge--Kutta
(DDISRK) methods of weak order one and two are calculated. Their
asymptotic stability as well as mean--square stability
(MS--stability) properties are studied for a linear stochastic
test equation with multiplicative noise. The stability functions
for the DDISRK methods are determined and their domains of
stability are compared to the corresponding domain of stability of
the considered test equation. Stability regions are presented for
various coefficients of the families of DDISRK methods in order to
determine step size restrictions such that the numerical
approximation reproduces the characteristics of the solution
process.
\end{abstract}
\end{frontmatter}
\section{Introduction} \label{Sec:Intro}
Numerical methods are an important tool for the calculation of
approximate solutions of stochastic differential equations (SDEs)
which possess no analytical solution formula. Therefore, many
approximation schemes have been developed in recent years and much
research has been carried out to develop derivative free
stochastic Runge--Kutta (SRK) type methods
\cite{BuTi04,Roe06a,Roe06b,Roe04,ToVA02}. Similar to the well
understood deterministic setting of ordinary differential
equations (ODEs), one has to pay much attention to the stability
properties of the solution as well as of the numerical
approximations. Therefore, implicit methods have been proposed for
the strong pathwise approximation of solutions of SDEs and their
stability has been analyzed \cite{BuTi04,Hig00,Pet98,Schu99}.
However, for the approximation of moments of the solution process
special numerical methods converging in the weak sense have to be
applied (see, e.g.,
\cite{KP99,Ko07,KoMi95,Mil04,Roe06a,Roe06b,Roe04,Roe06c,Roe07b,ToVA02})
and a stability analysis has to be carried out similar to that for
strong approximations \cite{Hig00,KoMi95,To05}. In the present
paper, we present families of first and second order diagonally
drift--implicit SRK (DDISRK) methods for the weak approximation of
SDEs contained in the class of SRK methods proposed by
R\"{o}{\ss}ler~\cite{Roe06c}. Further, we analyze their asymptotic
stability and mean--square stability for linear test equations
with multiplicative noise. Finally, the regions of stability of
the DDISRK methods are compared to the regions of stability of the
linear test equation. Thus, in Section~\ref{Sec:sec-ord-SRK} we
consider the class of SRK methods and coefficients families for
weak order one and order two DDISRK methods are presented. Then,
we discuss the concepts of stability for solutions of SDEs and for
numerical approximations in Section~\ref{Sec:Stab-ana-SDE} and
Section~\ref{Sec:Numer-mean-square-stab}, respectively. In
Section~\ref{Sec:Num-experiments}, some numerical experiments are
carried out in order to justify our theoretical results.

Let $(\Omega, \mathcal{F}, \Prob)$ be a probability space with a
filtration $(\mathcal{F}_t)_{t \geq 0}$ which fulfills the usual
conditions and $\mathcal{I}=[t_0,T]$ for some $0 \leq t_0 < T <
\infty$. Then, let $X=(X_t)_{t \in \mathcal{I}}$ denote the
solution process of an It\^{o} SDE
\begin{equation} \label{Ito-SDE-1}
    {\mathrm{d}} X_t = a(t,X_t) \, {\mathrm{d}}t + b(t,X_t) \,
    {\mathrm{d}}W_t, \qquad X_{t_0} = x_{0},
\end{equation}
for $t \in \mathcal{I}$ where $a : \mathcal{I} \times \mathbb{R}^d
\to \mathbb{R}^d$ is the drift and $b : \mathcal{I} \times
\mathbb{R}^d \to \mathbb{R}^{d \times m}$ is the diffusion,
$(W_t)_{t \geq 0}$ is a $m$-dimensional Wiener process and a
$\mathcal{F}_{t_0}$--measurable initial condition $x_0$
independent of $W_t-W_{t_0}$ for $t \geq t_0$ such that
$\E(\|x_0\|^{2r}) < \infty$ for some $r \in \mathbb{N}$. The $j$th
column of the $d \times m$--diffusion matrix $b=(b^{ij})$ will be
denoted by $b^j$ in the following. Further, we suppose that the
conditions of the existence and uniqueness theorem \cite{KP99} are
fulfilled for SDE~(\ref{Ito-SDE-1}).

In the following, we consider time discrete approximations
$Y^h=(Y_{t})_{t \in \mathcal{I}_h}$ w.r.t.\ a constant step size
$h=\frac{T-t_0}{N}$ for some $N \in \mathbb{N}$ and $\mathcal{I}_h
= \{t_0, t_1, \ldots, t_N\}$ where $t_n = t_0 + n \, h$ for $0
\leq n \leq N$. As usual, we also write $Y_n=Y_{t_n}$ for $0 \leq
n \leq N$. Further, let $C_P^l(\mathbb{R}^d, \mathbb{R})$ denote
the space of all $g \in C^l(\mathbb{R}^d,\mathbb{R})$ fulfilling a
polynomial growth condition \cite{KP99}.

\begin{Def}
A time discrete approximation $Y^h$ converges weakly with order
$p>0$ to $X$ as $h \rightarrow 0$ at time $t \in \mathcal{I}_h$ if
for each $f \in C_P^{2(p+1)}(\mathbb{R}^d, \mathbb{R})$ exists a
constant $C_f$, which does not depend on $h$, and a finite
$\delta_0 > 0$ such that for each $h \in \, ]0,\delta_0[\,$
\begin{equation}
    | \E(f(X_t)) - \E(f(Y_t)) | \leq C_f \, h^p \, .
\end{equation}
\end{Def}
\section{Diagonally Drift-Implicit Stochastic Runge--Kutta Methods}
\label{Sec:sec-ord-SRK}
For the weak approximation of the solution of the It\^{o}
SDE~(\ref{Ito-SDE-1}), we consider the class of SRK methods
introduced by R\"{o}{\ss}ler~\cite{Roe06c}. Then, the $d$-dimensional
approximation process $Y^h$ with $Y_n=Y_{t_n}$ for $t_n \in
\mathcal{I}_h$ is given by the following SRK method of $s$-stages
with $Y_0=x_0$ and
\begin{equation} \label{SRK-method-Ito-Wm-allg01}
    \begin{split}
    Y_{n+1} = Y_n & + \sum_{i=1}^s
    \alpha_i \, a(t_n+c_i^{(0)} h_n, H_i^{(0)}) \, h_n \\
    & + \sum_{i=1}^s
    \sum_{k=1}^m
    \big( {\beta_i^{(1)}} \, \hat{I}_{(k),n}
    + {\beta_i^{(2)}} \, \frac{\hat{I}_{(k,k),n}}{\sqrt{h_n}}
    \big) \, b^{k}(t_n+c_i^{(1)} h_n, H_i^{(k)}) \\
    & +
    \sum_{i=1}^s \sum_{k=1}^m
    \big( {\beta_i^{(3)}} \, \hat{I}_{(k),n}
    + {\beta_i^{(4)}} \, \sqrt{h_n} \big)
    \, b^{k}(t_n+{c}_i^{(2)} h_n, \hat{H}_i^{(k)})
    \end{split}
\end{equation}
for $n=0,1, \ldots, N-1$ with stage values
\begin{alignat*}{3}
    H_i^{(0)} &&= Y_n &+ \sum_{j=1}^{s} A_{ij}^{(0)}
    \, a(t_n+c_j^{(0)} h_n, H_j^{(0)}) \, h_n \\
    && &+ \sum_{j=1}^{s} \sum_{l=1}^m
    {B_{ij}^{(0)}} \, b^l(t_n+c_j^{(1)} h_n, H_j^{(l)}) \, \hat{I}_{(l),n} \\
    H_i^{(k)} &&= Y_n &+ \sum_{j=1}^{s} A_{ij}^{(1)}
    \, a(t_n+c_j^{(0)} h_n, H_j^{(0)}) \, h_n \\
    && &+ \sum_{j=1}^{s}
    {B_{ij}^{(1)}} \, b^k(t_n+c_j^{(1)} h_n, H_j^{(k)}) \,
    \sqrt{h_n} \\
    \hat{H}_i^{(k)} &&= Y_n &+ \sum_{j=1}^{s} {A}_{ij}^{(2)}
    \, a(t_n+c_j^{(0)} h_n, H_j^{(0)}) \, h_n \\
    && &+ \sum_{j=1}^{s} \sum_{\substack{l=1 \\ l \neq k}}^m
    {{B}_{ij}^{(2)}} \, b^l(t_n+c_j^{(1)} h_n, H_j^{(l)}) \,
    \frac{\hat{I}_{(k,l),n}}{\sqrt{h_n}}
\end{alignat*}
for $i=1, \ldots, s$ and $k=1, \ldots, m$. The random variables of
the method are defined by
\begin{equation}
    \hat{I}_{(k,l),n} = \begin{cases} \tfrac{1}{2} ( \hat{I}_{(k),n}
    \hat{I}_{(l),n} - \sqrt{h_n} \tilde{I}_{(k),n} ) & \text{if } k <
    l \\
    \tfrac{1}{2} ( \hat{I}_{(k),n}
    \hat{I}_{(l),n} + \sqrt{h_n} \tilde{I}_{(l),n} ) & \text{if } l < k \\
    \tfrac{1}{2} ( \hat{I}_{(k),n}^2 - h_n ) & \text{if } k =
    l
    \end{cases}
\end{equation}
for $1 \leq k,l \leq m$ with independent random variables
$\hat{I}_{(k),n}$ for $1 \leq k \leq m$ and $\tilde{I}_{(k),n}$
for $1 \leq k \leq m-1$ and $0 \leq n < N$. Thus, only $2m-1$
independent random variables have to be simulated for each step.
In the following, we choose $\hat{I}_{(k),n}$ as a three point
distributed random variable with $\Prob(\hat{I}_{(k),n} = \pm
\sqrt{3 \, h_n} ) = \frac{1}{6}$ and $\Prob(\hat{I}_{(k),n} = 0 )
= \frac{2}{3}$. The random variables $\tilde{I}_{(k),n}$ are
defined by a two point distribution with $\Prob(\tilde{I}_{(k),n}
= \pm \sqrt{h})=\tfrac{1}{2}$.

The main advantage of this class of SRK methods is the significant
reduction of complexity compared to present SRK methods in recent
literature, because the number of stages does not depend on the
dimension $m$ of the driving Wiener process \cite{Roe06c}.
We denote by $\alpha=(\alpha_i)$ and $\beta^{(k)}=(\beta^{(k)}_i)$
for $1 \leq k \leq 4$ the corresponding vectors of weights and by
$A^{(k)} = (A^{(k)}_{ij})$ and $B^{(k)}=(B^{(k)}_{ij})$ for
$k=0,1,2$ the corresponding coefficients matrices. Then, the
coefficients of the SRK method~(\ref{SRK-method-Ito-Wm-allg01}) can
be represented by an extended Butcher array:
\begin{equation*}
{
\begin{tabular}{c|c|c|c}
    $c^{(0)}$ & ${A}^{(0)}$ & ${B^{(0)}}$ & \\
    \cline{1-4}
    $c^{(1)}$ & $A^{(1)}$ & ${B^{(1)}}$ & \\
    \cline{1-4}
    ${c}^{(2)}$ & $A^{(2)}$ & $B^{(2)}$ & \\
    \hline
    & $\alpha^T$ & ${\beta^{(1)}}^T$ & ${\beta^{(2)}}^T$ \\
    \cline{2-4}
    & & ${\beta^{(3)}}^T$ & ${\beta^{(4)}}^T$
\end{tabular}
}
\end{equation*}
Weak order one and two conditions for the SRK method
(\ref{SRK-method-Ito-Wm-allg01}) have been calculated in
\cite{Roe06c} by applying the colored rooted tree theory due to
R\"{o}{\ss}ler introduced in~\cite{Roe06a}.
Now, let $p_S=p$ denote the order of convergence of the SRK method
(\ref{SRK-method-Ito-Wm-allg01}) if it is applied to an SDE and let
$p_D$ with $p_D \geq p_S$ denote the order of convergence if it is
applied to a deterministic ODE, i.e., SDE~(\ref{Ito-SDE-1}) with $b
\equiv 0$ and we also write $(p_D,p_S)$ \cite{Roe06a,Roe06b}.
Since we are interested in SRK methods which inherit good
stability properties, we consider families of DDISRK methods which
are diagonally implicit in the deterministic part of the scheme.
\subsection{Weak Order One DDISRK Methods}
Firstly, we consider weak order one DDISRK methods
(\ref{SRK-method-Ito-Wm-allg01}) with $s=1$ stage \cite{Roe06c}.
However, in order to cover the stochastic theta method
\cite{Hig00,KP99}, we also consider the case that the stage number
is $s=2$ for the drift function only, whereas it is still one for
the diffusion function. Then, from the order conditions
\cite{Roe06c} it follows that the family of weak order one DDISRK
methods is characterized by the Butcher table
(\ref{SRK-DDI-Coefficients-order1})
\begin{table}
\begin{equation} \label{SRK-DDI-Coefficients-order1}
{
\begin{tabular}{c|cc|cc|cc}
    $c_1$ & $c_1$ & 0 & 0 & 0 & \\
    $c_2+c_3$ & $c_2$ & $c_3$ & $c_4$ & 0 & \\
    \cline{1-5}
    0 & 0 & 0 & 0 & 0 & \\
    0 & 0 & 0 & 0 & 0 & \\
    \cline{1-5}
    0 & 0 & 0 & 0 & 0 & \\
    0 & 0 & 0 & 0 & 0 & \\
    \hline
    & $1-c_5$ & $c_5$ & $1$ & 0 & 0 & 0 \\
    \cline{2-7}
    & & & 0 & 0 & 0 & 0
\end{tabular}
}
\end{equation}
\end{table}
with some coefficients $c_1, c_2, c_3, c_4, c_5 \in \mathbb{R}$.
As an example, in the case of $s=1$ stage we obtain for
$c_1=c_2=c_3=c_4=c_5=0$ the explicit Euler-Maruyama scheme of
order $(1,1)$ \cite{KP99}. For $c_1=\frac{1}{2}$ and
$c_2=c_3=c_4=c_5=0$ we obtain the SRK scheme DDIRDI1 with $s=1$
stage of order $(2,1)$, which reduces to the midpoint rule if it
is applied to an ODE \cite{HW96}. If we consider the case of $s=2$
stages, then we get for $c_1=0$, $c_2=1-\theta$, $c_3=\theta$,
$c_4=1$ and $c_5=\theta$ for some $\theta \in [0,1]$ the SRK
scheme DDIRDI2 of order $(2,1)$ which coincides with the
stochastic theta method \cite{Hig00,KP99,SaMi96,Schu96}. Further,
for $s=2$ stages with $c_1=c_3=\frac{3+\sqrt{3}}{6}$,
$c_2=-\frac{1}{\sqrt{3}}$, $c_5=\frac{1}{2}$ and $c_4 \in
\mathbb{R}$ we get DDISRK schemes of order $(3,1)$ which are
A-stable in case of ODEs \cite{HW96}. Especially, in the case of
$c_4=\frac{3}{2}$ we denote the scheme as DDIRDI3 in the
following. Note that the schemes DDIRDI1 and DDIRDI2 for
$\theta=\frac{1}{2}$ are also A-stable if they are applied to ODEs
\cite{HW96}.
\subsection{Weak Order Two DDISRK Methods}
Next, we consider weak second order DDISRK methods
(\ref{SRK-method-Ito-Wm-allg01}) with $s=3$ stages. Here, we claim
that $\alpha_3=0$ and $A^{(2)}_{ij}=0$ for $1 \leq i,j \leq 3$ in
order to reduce the computational effort. Since in this case the
third stage $H^{(0)}_3$ does not matter anymore, we let
$A^{(0)}_{3j}=B^{(0)}_{3j}=0$ for $1 \leq j \leq 3$.
Then, we can obtain $\alpha_1=\alpha_2=\frac{1}{2}$ from some
order conditions of weak order two \cite{Roe06c} and from the
classification given in \cite{DeRoe06d} in the case of an explicit
SRK scheme. On the other hand, since we assume $A^{(2)} \equiv 0$,
all conditions for $A^{(0)}$ to satisfy are $\alpha^T A^{(0)} e =
\frac{1}{2}$ only. Therefore, we can consider arbitrary
coefficients $A^{(0)}_{ij}$ as long as $\alpha^T A^{(0)} e =
\frac{1}{2}$ is fulfilled for $\alpha_1=\alpha_2=\frac{1}{2}$.
As a result of this, the weak order two DDISRK schemes
(\ref{SRK-method-Ito-Wm-allg01}) are given by the infinite
coefficients family (\ref{SRK-DDI-Coefficients})
\begin{table}
\begin{equation} \label{SRK-DDI-Coefficients}
{
\begin{tabular}{c|ccc|ccc|ccc}
    $c_1$ & $c_1$ & 0 & 0 & 0 & 0 & 0 & \\
    $1-c_1$ & $1-c_1-c_2$ & $c_2$ & 0 & 1 & 0 & 0 & \\
    0 & 0 & 0 & 0 & 0 & 0 & 0 \\
    \cline{1-7}
    0 & 0 & 0 & 0 & 0 & 0 & 0 & \\
    $c_3^2$ & $c_3^2$ & 0 & 0 & $c_3$ & 0 & 0 & \\
    $c_3^2$ & $c_3^2$ & 0 & 0 & $-c_3$ & 0 & 0 & \\
    \cline{1-7}
    0 & 0 & 0 & 0 & 0 & 0 & 0 & \\
    0 & 0 & 0 & 0 & $c_4$ & 0 & 0 & \\
    0 & 0 & 0 & 0 & $-c_4$ & 0 & 0 & \\
    \hline
    & $\frac{1}{2}$ & $\frac{1}{2}$ & 0 & $1-\frac{1}{2c_3^2}$ &
    $\frac{1}{4c_3^2}$ & $\frac{1}{4c_3^2}$ & 0 & $\frac{1}{2c_3}$ &
    $-\frac{1}{2c_3}$ \\
    \cline{2-10}
    & & & & $-\frac{1}{2c_4^2}$ & $\frac{1}{4c_4^2}$ &
    $\frac{1}{4c_4^2}$ & 0 & $\frac{1}{2c_4}$ & $-\frac{1}{2c_4}$
\end{tabular}
}
\end{equation}
\end{table}
with $c_1,c_2 \in \mathbb{R}$ and $c_3,c_4 \in \mathbb{R}
\setminus \{0\}$. Clearly, one has to solve 2 (in general
nonlinear) systems of equations from the stage values, each of
dimension $d$, for the DDISRK method if $c_1 \neq 0$ and $c_2 \neq
0$. Therefore, as in the deterministic setting, some simplified
Newton iterations have to be performed in each step in order to
solve the nonlinear system of equations \cite{DeKv07,HW96}. As an
example, for $c_1=c_2=\frac{3+\sqrt{3}}{6}$ and for all $c_3, c_4
\in \mathbb{R} \setminus \{0\}$ we obtain an DDISRK scheme of
order $(3,2)$ which is A-stable if it is applied to a
deterministic ODE \cite{HW96}.
\section{Stability analysis for SDEs}
\label{Sec:Stab-ana-SDE}
For SDEs several stochastic stability concepts have been proposed
in literature, see e.g.,
\cite{Hig00,HoPl94,KP99,KoMi95,MilPlSchu98,SaMi93,SaMi96,Schu96,To05}
and the literature therein.
In the following, we consider SDE (\ref{Ito-SDE-1}) with a steady
solution $X_t \equiv 0$ such that $a(t,0)=b(t,0)=0$ holds, which
is also called an equilibrium position. Suppose that there exists
a unique solution $X_t = X(t;t_0,x_0)$ for all $t \geq t_0$ and
for each nonrandom initial value $x_0$ under consideration. Then,
stochastic stability can be defined as the stochastic counterparts
of stability, asymptotic stability and asymptotic stability in the
large for ODEs \cite{BuBuMi00,Has80,KP99}.
\begin{Def} \label{Def:stoch-stability-sol}
    Let $X$ be the solution of the scalar It\^{o} SDE (\ref{Ito-SDE-1}).
    Then, the equilibrium position of the SDE is said to be
    \begin{enumerate}[(i)]
        \item \label{Def:stoch-stability-sol:a}
        stochastically stable if for all $\epsilon > 0$
        and $t_0 \geq 0$ holds
        \[
            \lim_{x_0 \to 0} \Prob \left( \sup_{t \geq t_0}
            |X(t;t_0,x_0)| \geq \epsilon \right) = 0 \, ,
        \]
        \item \label{Def:stoch-stability-sol:b}
        stochastically asymptotically stable if
        (\ref{Def:stoch-stability-sol:a}) holds and if
        \[
            \lim_{x_0 \to 0} \Prob \left( \lim_{t \to \infty}
            |X(t;t_0,x_0)| =0 \right) = 1 \, ,
        \]
        \item \label{Def:stoch-stability-sol:c}
        or stochastically asymptotically stable in the large if
        (\ref{Def:stoch-stability-sol:a}) holds and if for all $x_0 \in
        \mathbb{R}$ holds
        \[
            \Prob \left( \lim_{t \to \infty} |X(t;t_0,x_0)| = 0
            \right) = 1 \, .
        \]
    \end{enumerate}
\end{Def}
Further, stability analysis involving the $p$th moments of the
solution process is also widely considered, see e.g.\
\cite{BuBuMi00,Hig00,HoPl94,KP99,KoMi95}.
\begin{Def} \label{Def:ms-stability-sol}
    Let $X$ be the solution of the scalar It\^{o} SDE (\ref{Ito-SDE-1}).
    Then, the equilibrium position of the SDE is said to be
    \begin{enumerate}[(i)]
        \item \label{Def:ms-stability-sol:a}
        stable in the $p$th--mean if for every $\epsilon>0$ and
        $t_0 \geq 0$ there exists a $\delta=\delta(t_0,\epsilon) >
        0$ such that for all $t \geq t_0$ and $|x_0| < \delta$
        \[
            \E \left( |X(t;t_0,x_0)|^p \right) < \epsilon \, ,
        \]
        \item \label{Def:ms-stability-sol:b}
        asymptotically stable in the $p$th--mean if
        (\ref{Def:ms-stability-sol:a}) holds and if there exists a
        $\delta_0=\delta_0(t_0)>0$ such that for all $|x_0| < \delta_0$
        \[
            \lim_{t \to \infty} \E \left( |X(t;t_0,x_0)|^p \right)
            = 0 \, .
        \]
    \end{enumerate}
\end{Def}
The most frequently used cases in applications are $p=1$ and $p=2$,
i.e., stability in mean (M-stability) and mean--square stability
(MS-stability). In the present paper, we will focus on asymptotic
stability in the large and MS-stability for a linear test equation
with multiplicative noise \cite{HeSp93,Hig00,HoPl94,SaMi96}
\begin{equation} \label{Lin-stoch-test-eqn}
    {\mathrm{d}} X_t = \lambda \, X_t \, {\mathrm{d}}t + \mu \,
    X_t \, {\mathrm{d}}W_t
\end{equation}
for $t \geq t_0$ and with some constants $\lambda,\mu \in
\mathbb{C}$ and with a nonrandom initial condition $X_{t_0}=x_0 \in
\mathbb{R} \setminus \{0\}$, which reproduces the dynamics of more
complex SDEs better than in the case of additive noise
\cite{HeSp92,KP99}.
The exact solution of (\ref{Lin-stoch-test-eqn}) can be calculated
as $X_t = x_0 \, \exp ((\lambda-\tfrac{1}{2} \mu^2) (t-t_0) + \mu
(W_t-W_{t_0}))$ which is stochastically asymptotically stable in
the large \cite{SaMi96} if
\begin{equation}
    \begin{split}
        \lim_{t \to \infty} |X_t|=0 \text{ with
        probability 1} \quad &\Leftrightarrow \quad \Re(\lambda-\tfrac{1}{2}
        \mu^2) < 0 \, .
    \end{split}
\end{equation}
We calculate that $|X_t|^p = |x_0|^p \, \exp ( p \, \Re (\lambda -
\tfrac{1}{2} \mu^2 ) (t-t_0) + p \, \Re(\mu) ( W_t - W_{t_0} ))$
which yields $\E ( |X_t|^p ) = |x_0|^p \, \exp ( p \, \Re (
\lambda - \tfrac{1}{2} \mu^2 ) (t-t_0) + \tfrac{1}{2} p^2
(\Re(\mu) )^2 (t-t_0) )$. Then the $p$th--mean stability domain
where SDE (\ref{Lin-stoch-test-eqn}) possesses an equilibrium
position can be determined as follows:
\begin{equation} \label{Stability-domain-pth-mean}
    \lim_{t \to \infty} \E ( |X_t|^p ) = 0 \quad
    \Leftrightarrow \quad 2 \, \Re(\lambda) - \Re(\mu^2) + p \, (\Re(\mu))^2 < 0 \, .
\end{equation}
Thus, the equilibrium position of SDE (\ref{Lin-stoch-test-eqn}) is
asymptotically MS--stable if
\begin{equation} \label{Stability-domain-MS-mean-Loes}
    \lim_{t \to \infty} \E ( |X_t|^2 ) = 0 \quad
    \Leftrightarrow \quad 2 \, \Re(\lambda) + |\mu|^2 < 0
\end{equation}
for the coefficients $\lambda, \mu \in \mathbb{C}$ (see, e.g.,
\cite{Hig00,SaMi96,To05}). We remark that due to $\Re(2 \lambda -
\mu^2) \leq 2 \, \Re(\lambda) + |\mu|^2$ MS-stability always induces
asymptotically stability in the large. Further, for $\mu=0$ the
stability condition (\ref{Stability-domain-pth-mean}) reduces to the
well known deterministic stability condition $\Re(\lambda)<0$.
\section{Numerical stability of SRK methods}
\label{Sec:Numer-mean-square-stab}
We are now looking for conditions such that a numerical method
applied to SDE (\ref{Lin-stoch-test-eqn}) yields numerically
stable solutions. Therefore, we say that the method is numerically
asymptotically stable or MS--stable if the numerical solutions
$Y_n$ satisfy $\lim_{n \to \infty} |Y_n|=0$ with probability one
or $\lim_{n \to \infty} \E \left( |Y_{n}|^2 \right) = 0$,
respectively. If we apply the numerical method to the linear test
equation (\ref{Lin-stoch-test-eqn}), then we obtain with the
parametrization $\hat{h} = \lambda \, h$ and $k = \mu \sqrt{h}$
\cite{Hig00,KoMi95} a one--step difference equation of the form
\begin{equation} \label{Lemma-recursion-one-step-method}
    Y_{n+1} = R_n(\hat{h}, k) \, Y_n
    = \prod_{i=0}^{n} R_i(\hat{h}, k) \, Y_0 \, .
\end{equation}
with a stability function $R_n(\hat{h}, k)$. The domain of
asymptotic stability of a numerical method can be determined by
the following lemma \cite{Hig00}:
\begin{Lem} \label{Lemma-asym-stab}
    Given a sequence of real-valued, non-negative, independent and
    identically distributed random variables $(|R_n(\hat{h}, k)|)_{n \in
    \mathbb{N}_0}$, consider the sequence of random variables
    $(|Y_n|)_{n \in \mathbb{N}_0}$ defined by
    (\ref{Lemma-recursion-one-step-method})
    where $|Y_0| \neq 0$ with probability 1.
    Suppose that the random variables $\log(R_n(\hat{h}, k))$ are
    square-integrable. Then
    \begin{equation}
        \lim_{n \to \infty} |Y_n| =0 \, , \text{ with
        probability 1 } \Leftrightarrow \, \E(\log(R_n(\hat{h}, k))) < 0 \, .
    \end{equation}
\end{Lem}
We call the set $\mathcal{R}_{AS} = \{ (\hat{h}, k) \in
\mathbb{C}^2 : \E(\log(R_n(\hat{h}, k))) < 0 \} \subset
\mathbb{C}^2$ the domain of asymptotical stability of the method.
Note that one can also find some alternative parameterizations
like $k = -\frac{\mu^2}{\lambda}$ in the literature
\cite{BuBuMi00,SaMi96,To05}.
Analogously, if we calculate the mean--square norm $z_n =
\E(|Y_{n}|^2)$ then we obtain a one--step difference equation of
the form $z_{n+1} = \hat{R}(\hat{h}, k) \, z_n$ where
$\hat{R}(\hat{h},k)=\E(|R_n(\hat{h}, k)|^2)$ is called the
MS--stability function of the numerical method. Thus, we obviously
yield MS--stability, i.e. $z_n \to 0$ as $n \to \infty$, if
$\hat{R}(\hat{h},k)<1$. The set $\mathcal{R}_{MS} = \{ (\hat{h},k)
\in \mathbb{C}^2 : \hat{R}(\hat{h},k)<1 \} \subset \mathbb{C}^2$
is called the domain of MS--stability of the method.

Especially, the domain is called region of stability in the case
of $(\hat{h},k) \in \mathbb{R}^2$ \cite{SaMi96}. The numerical
method is said to be $A$--stable if the domain of stability of the
test equation (\ref{Lin-stoch-test-eqn}) is a subset of the domain
of numerical stability.
Since the domain of stability for $\lambda, \mu \in \mathbb{C}$ is
not easy to visualize, we restrict our attention to figures
presenting the region of stability with $\lambda, \mu \in
\mathbb{R}$ in the $\hat{h}$--$k^2$ plane. Then, for fixed values
of $\lambda$ and $\mu$, the set $\{(\lambda \, h, \mu^2 \, h)
\subset \mathbb{R}^2 : h > 0\}$ is a straight ray starting at the
origin and going through the point $(\lambda, \mu^2)$. Clearly,
varying the step size $h$ corresponds to moving along this ray.
For $\lambda,\mu \in \mathbb{R}$, the region of asymptotical
stability for SDE (\ref{Lin-stoch-test-eqn}) reduces to the area
of the $\hat{h}$--$k^2$ plane with the $\hat{h}$--axis as the
lower bound if $\hat{h} < 0$ and with $k^2>2 \hat{h}$ as the lower
bound if $\hat{h} \geq 0$ whereas the region of MS--stability for
SDE (\ref{Lin-stoch-test-eqn}) reduces to the area of the
$\hat{h}$--$k^2$ plane with the $\hat{h}$--axis as the lower bound
and $k^2<-2 \hat{h}$ as the upper bound for $\hat{h} < 0$.

In the following, all figures presenting regions of stability for
some numerical method under consideration are plotted by the
software Mathematica. The regions of numerical asymptotically
stability and MS--stability are indicated by two dark--grey tones
whereas the regions of MS--stability are more dark than the
regions of asymptotical stability. Further, the corresponding
regions of stability for the test equation
(\ref{Lin-stoch-test-eqn}) are filled by two light--grey tones
whereas again the regions of MS--stability are more dark than the
regions of asymptotical stability. In all presented figures, the
regions of MS--stability are a subset of the regions of
asymptotical stability.
\subsection{Stability of Order One DDISRK Schemes}
\label{Sec:MS-stab-oder-one-SRK}
\begin{figure}[tbp]
\begin{center}
\includegraphics[width=6.8cm]{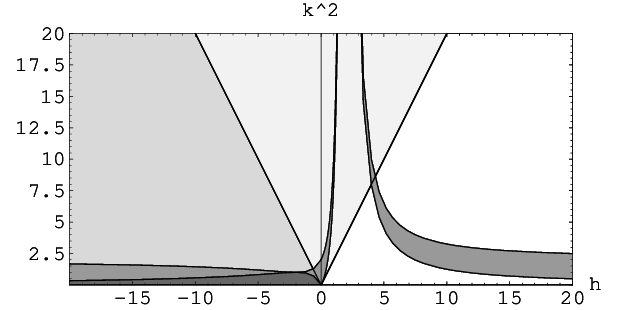}
\includegraphics[width=6.8cm]{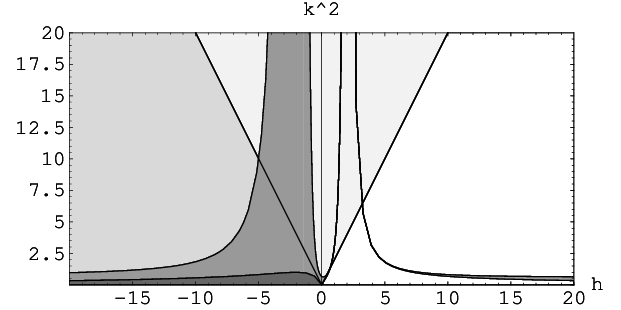}
\caption{Asymptotical and mean--square stability region for
DDIRDI1 with two point and three point distributed random
variables in the left and right figure, respectively.}
\label{Bild-DDIRDI1}
\end{center}
\end{figure}

\begin{figure}[tbp]
\begin{center}
\includegraphics[width=6.8cm]{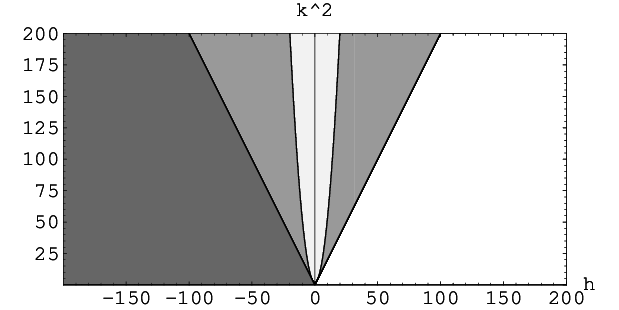}
\includegraphics[width=6.8cm]{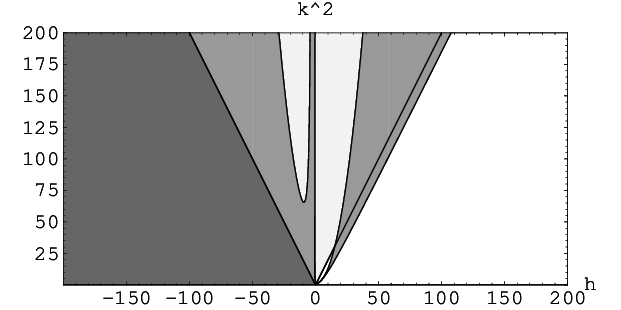}
\caption{Asymptotical and mean--square stability region for
DDIRDI2 with two point and three point distributed random
variables in the left and right figure, respectively.}
\label{Bild-DDIRDI2}
\end{center}
\end{figure}

\begin{figure}[tbp]
\begin{center}
\includegraphics[width=6.8cm]{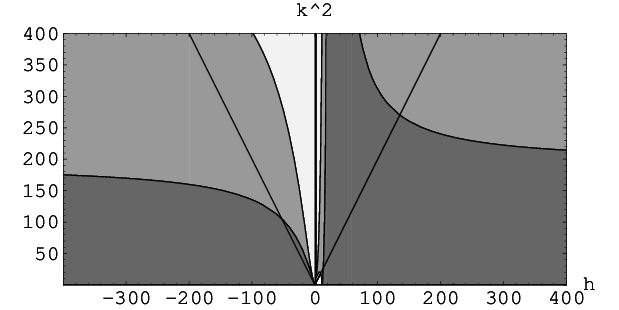}
\includegraphics[width=6.8cm]{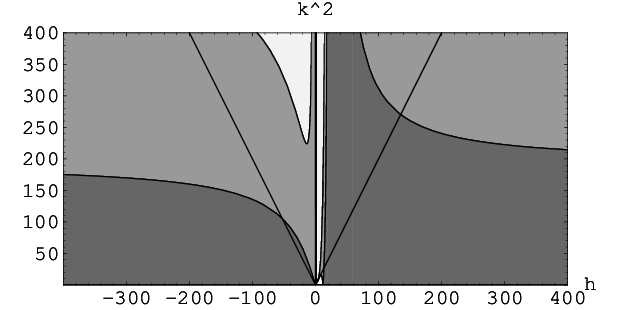}
\caption{Asymptotical and mean--square stability region for
DDIRDI3 with two point and three point distributed random
variables in the left and right figure, respectively.}
\label{Bild-DDIRDI3}
\end{center}
\end{figure}

We consider the family of order one DDISRK schemes
(\ref{SRK-DDI-Coefficients-order1}) with coefficients $c_1,
\ldots, c_5 \in \mathbb{R}$. If we apply these schemes to the
linear test equation (\ref{Lin-stoch-test-eqn}) then we obtain the
difference equation
\begin{equation} \label{SRK-order1-test-eqn}
    Y_{n+1} = Y_n + (1-c_5) \, \lambda \, h \, H_1^{(0)}
    + c_5 \, \lambda \, h \, H_2^{(0)} + \mu \, \hat{I}_{(1),n} \,
    H^{(1)}_1
\end{equation}
with the stage values
\begin{equation}
    \begin{split}
    H^{(0)}_1 &= Y_n + c_1 \, \lambda \, h \, H^{(0)}_1 \\
    H^{(0)}_2 &= Y_n + c_2 \, \lambda \, h \, H^{(0)}_1 +
    c_3 \, \lambda \, h \, H^{(0)}_2 + c_4 \, \mu \, \hat{I}_{(1),n} \,
    H^{(1)}_1 \\
    H^{(1)}_1 &= Y_n
    \end{split}
\end{equation}
where the implicit equations for $H^{(0)}_1$ and $H^{(0)}_2$ can
be solved in the case of $1-c_1 \, \lambda \, h \neq 0$ and $1-c_3
\, \lambda \, h \neq 0$, which is fulfilled for step sizes $h \neq
\frac{1}{c_1 \, \lambda}$ if $c_1 \neq 0$ and $h \neq \frac{1}{c_3
\, \lambda}$ if $c_3 \neq 0$. With $\hat{h} = \lambda \, h$ and $k
= \mu \sqrt{h}$ let
\begin{equation*}
    \begin{split}
    &\Gamma = 1 + \frac{\hat{h}-c_3 \, \hat{h}^2
    + c_5 \, (c_2+c_3-c_1) \, \hat{h}^2}{(1-c_1
    \, \hat{h})(1-c_3 \, \hat{h})}
    \, ,
    \qquad
    \Sigma = \frac{c_4 \, c_5 \, \hat{h}}{1-c_3 \, \hat{h}}
    \, k + k \, .
    \end{split}
\end{equation*}
Then, we can write (\ref{SRK-order1-test-eqn}) by the recursion
formula $Y_{n+1} = R_n(\hat{h}, k) \, Y_n$ with the stability
function $R_n(\hat{h}, k) = \Gamma + h^{-1/2} \, \Sigma \,
\hat{I}_{(1),n}$ for $n=0, \ldots, N-1$. Since the SRK schemes
(\ref{SRK-DDI-Coefficients-order1}) are of weak order one, we can
substitute the tree point distributed random variables
$\hat{I}_{(j),n}$ by two point distributed random variables
$\tilde{I}_{(j),n}$ for $1 \leq j \leq m$ in
(\ref{SRK-method-Ito-Wm-allg01}) and consider $R_n(\hat{h}, k) =
\Gamma + h^{-1/2} \, \Sigma \, \tilde{I}_{(1),n}$ instead.
Now, we analyse the asymptotic stability of the SRK schemes
(\ref{SRK-DDI-Coefficients-order1}) by applying
Lemma~\ref{Lemma-asym-stab}. Further, in order to analyse the
MS-stability, we calculate the mean--square norm $z_n =
\E(|Y_{n}|^2)$. Then, we obtain the recursion formula $z_{n+1} =
\hat{R}(\hat{h}, k) \, z_n$ with the MS--stability function
$\hat{R}(\hat{h}, k) = |\Gamma|^2 + |\Sigma|^2$.
\begin{Sat} \label{Sat-1}
    For SDE (\ref{Lin-stoch-test-eqn}) with  $\lambda, \mu \in
    \mathbb{C}$, the SRK schemes (\ref{SRK-DDI-Coefficients-order1}) are
    \begin{enumerate}[(i)]
    \item numerically asymptotical stable if $|\Gamma^2 - 3 \, \Sigma^2 | \,
        |\Gamma|^4 < 1$
    and in the case that the random variables $\hat{I}_{(j),n}$ are replaced by
    $\tilde{I}_{(j),n}$ for $1 \leq j \leq m$,
    if $|\Gamma^2 - \Sigma^2 | < 1$,
    \item numerically MS--stable if $|\Gamma|^2 + |\Sigma|^2 < 1$.
    \end{enumerate}
\end{Sat}
Here, we have to point out that the distribution of the random
variables used for the numerical method has significant influence
on the domain of asymptotical stability.
As an example, for DDIRDI1 we have $\Gamma = \frac{1+\frac{1}{2}
\hat{h}}{1-\frac{1}{2} \hat{h}}$, $\Sigma = k$ and we calculate
$\mathcal{R}_{AS} = \{ (\hat{h}, k) \in \mathbb{C}^2 :
|(1+\frac{1}{2} \hat{h})^2 - (1-\frac{1}{2} \hat{h})^2 k^2| <
|1-\frac{1}{2} \hat{h}|^2 \}$ if the random variables
$\tilde{I}_{(j),n}$ are used and $\mathcal{R}_{AS} = \{ (\hat{h},
k) \in \mathbb{C}^2 : |(1-\frac{1}{2} \hat{h})^2 (1+\frac{1}{2}
\hat{h})^2 (1 - 3 k^2)| |1+\frac{1}{2} \hat{h}|^4 < |1-\frac{1}{2}
\hat{h}|^4 \}$ if $\hat{I}_{(j),n}$ are used. In both cases, we
get $\mathcal{R}_{MS} = \{ (\hat{h}, k) \in \mathbb{C}^2 : 2
\Re(\hat{h}) + ((1-\frac{1}{2} \Re(\hat{h}))^2 + \frac{1}{4}
(\Im(\hat{h}))^2) \, |k|^2 < 0\}$. The corresponding regions of
stability are presented for both cases in
Figure~\ref{Bild-DDIRDI1}. For DDIRDI2 we get $\Gamma =
\frac{1+(1-\theta) \hat{h}}{1-\theta \hat{h}}$ and $\Sigma =
\frac{\theta \hat{h}}{1-\theta \hat{h}} k + k$, $\theta \in
[0,1]$. Then, for $\theta=\frac{1}{2}$ follows $\mathcal{R}_{AS} =
\{ (\hat{h}, k) \in \mathbb{C}^2 : |(1+\frac{1}{2} \hat{h})^2 -
k^2| < |1-\frac{1}{2} \hat{h}|^2 \}$ if $\tilde{I}_{(j),n}$ are
used, $\mathcal{R}_{AS} = \{ (\hat{h}, k) \in \mathbb{C}^2 :
|(1-\frac{1}{2} \hat{h})^2 ((1+\frac{1}{2} \hat{h})^2 - 3 k^2)|
|1+\frac{1}{2} \hat{h}|^4 < |1-\frac{1}{2} \hat{h}|^4 \}$ if
$\hat{I}_{(j),n}$ are used and $\mathcal{R}_{MS} = \{ (\hat{h}, k)
\in \mathbb{C}^2 : 2 \Re(\hat{h}) + |k|^2 < 0\}$.
Thus, the scheme DDIRDI2 with $\theta=\frac{1}{2}$ is $A$--stable
w.r.t.\ MS--stability and the corresponding regions are presented
in Figure~\ref{Bild-DDIRDI2} (see also \cite{Hig00,SaMi96}).
Analogously, we can calculate the domains of stability for DDIRDI3
which are presented in Figure~\ref{Bild-DDIRDI3}. For all
considered schemes, we can see the influence of the random
variables used by the scheme to the domain of asymptotical
stability.
\subsection{Stability of Order Two DDISRK Schemes}
\label{Sec:MS-stab-oder-two-SRK}
Next, we apply the DDISRK method (\ref{SRK-method-Ito-Wm-allg01})
with the coefficients (\ref{SRK-DDI-Coefficients}) to the linear
test equation (\ref{Lin-stoch-test-eqn}). Then we obtain the
difference equation
\begin{equation} \label{SRK-order2-test-eqn}
    \begin{split}
    Y_{n+1} = \,\, &Y_n + \frac{1}{2} \, \lambda \, h \, H_1^{(0)} +
    \frac{1}{2} \, \lambda \, h \, H_2^{(0)} \\
    &+ \left( 1-\frac{1}{2c_3^2} \right) \, \mu \, \hat{I}_{(1)}
    \, H_1^{(1)} + \frac{1}{4 c_3^2} \, \mu \, \hat{I}_{(1)} \,
    H_2^{(1)} + \frac{1}{4 c_3^2} \, \mu \, \hat{I}_{(1)} \,
    H_3^{(1)} \\
    &+ \frac{1}{2c_3} \, \mu \, \frac{\hat{I}_{(1,1)}}{\sqrt{h}}
    \, H_2^{(1)} - \frac{1}{2c_3} \, \mu \, \frac{\hat{I}_{(1,1)}}{\sqrt{h}}
    \, H_3^{(1)}
    \end{split}
\end{equation}
with stage values
\begin{equation}
    \begin{split}
    H_1^{(0)} &= Y_n + c_1 \, \lambda \, h \, H_1^{(0)} \\
    H_2^{(0)} &= Y_n + (1-c_1-c_2) \, \lambda \, h \, H_1^{(0)} +
    c_2 \, \lambda \, h \, H_2^{(0)} + \mu \, \hat{I}_{(1)} \,
    H_1^{(1)} \\
    H_1^{(1)} &= Y_n \\
    H_2^{(1)} &= Y_n + c_3^2 \, \lambda \, h \, H_1^{(0)} + c_3 \,
    \mu \, \sqrt{h} \, H_1^{(1)} \\
    H_3^{(1)} &= Y_n + c_3^2 \, \lambda \, h \, H_1^{(0)} - c_3 \,
    \mu \, \sqrt{h} \, H_1^{(1)}
    \end{split}
\end{equation}
where the values $\hat{H}^{(1)}_i$ do not appear due to $m=1$ and
$A^{(2)} \equiv 0$.
Suppose that $1-c_1 \, \lambda \, h \neq 0$ and that $1-c_2 \,
\lambda \, h \neq 0$ which can always be fulfilled for step sizes
$h$ with $h \neq \frac{1}{c_1 \, \lambda}$ and $h \neq
\frac{1}{c_2 \, \lambda}$. Then the implicit equations for
$H_1^{(0)}$ and $H_2^{(0)}$ can always be solved. Let with
$\hat{h} = \lambda \, h$ and $k = \mu \sqrt{h}$
\begin{equation*}
    \begin{split}
    \Gamma &= 1 + \frac{\hat{h} + (\frac{1}{2}-c_1-c_2) \hat{h}^2}
    {(1-c_1 \, \hat{h})(1-c_2 \, \hat{h})} \, , \quad
    \Sigma = \frac{\hat{h} - \frac{1}{2}(c_1+c_2) \hat{h}^2}{(1-c_1 \,
    \hat{h})(1-c_2 \, \hat{h})} \, k + k \, , \quad
    \Lambda = \frac{1}{2} k^2 \, .
    \end{split}
\end{equation*}
Then, we yield for (\ref{SRK-order2-test-eqn}) the recursion
formula $Y_{n+1} = R_n(\hat{h},k) \, Y_n$ with the stability
function $R_n(\hat{h},k) = \Gamma - \Lambda + h^{-1/2} \, \Sigma
\, \hat{I}_{(1),n} + h^{-1} \, \Lambda \, \hat{I}_{(1),n}^2$.
In order to analyse the asymptotic stability of the weak order two
DDISRK schemes (\ref{SRK-DDI-Coefficients}) we apply again
Lemma~\ref{Lemma-asym-stab}. For the determination of the domain
of MS--stability, we calculate the mean--square norm of
(\ref{SRK-order2-test-eqn}). Then, we obatin the recursion
$z_{n+1} = \hat{R}(\hat{h}, k) \, z_n$ with MS--stability function
$\hat{R}(\hat{h}, k) = |\Gamma|^2+|\Sigma|^2+2 |\Lambda|^2$. Both
stability functions $R_n(\hat{h}, k)$ and $\hat{R}(\hat{h}, k)$
depend only on the coefficients $c_1$ and $c_2$ of the scheme,
i.e. the coefficients $c_3$ and $c_4$ are not relevant for the
stability in the case of the scalar linear test equation
(\ref{Lin-stoch-test-eqn}).
\begin{Sat} \label{Sat-Order2-Stability-DDISRK}
    For SDE (\ref{Lin-stoch-test-eqn}) with  $\lambda, \mu \in
    \mathbb{C}$, the SRK schemes (\ref{SRK-DDI-Coefficients}) are
    \begin{enumerate}[(i)]
    \item numerically asymptotical stable if $|(\Gamma+2 \,
    \Lambda)^2 - 3 \, \Sigma^2 | \, |\Gamma-\Lambda|^4 < 1$,
    \item numerically MS--stable if $|\Gamma|^2+|\Sigma|^2+2
    |\Lambda|^2 < 1$.
    \end{enumerate}
\end{Sat}
\begin{figure}[tbp]
\begin{center}
\includegraphics[width=6.8cm]{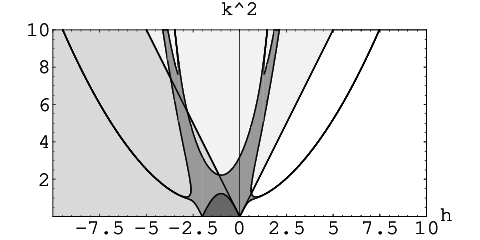}
\includegraphics[width=6.8cm]{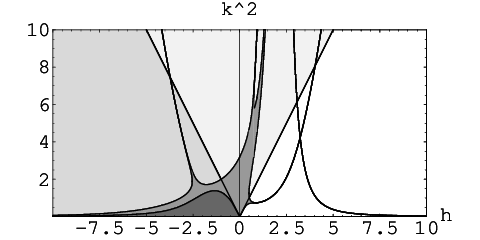}
\caption{Asymptotical and mean--square stability regions for RI6
with $c_1 = c_2 = 0$ on the left and for DDIRDI4 with $c_1=0$ and
$c_2=\frac{1}{2}$ on the right.} \label{Bild-RI6-DDIRDI4}
\end{center}
\end{figure}
If we choose $c_1=c_2=0$ and $c_3=c_4=1$ in
(\ref{SRK-DDI-Coefficients}), then we yield the explicit SRK
scheme RI6 calculated in \cite{Roe06c} which coincides for $m=1$
with the SRK scheme due to Platen \cite{KP99,Roe06b,To05}. If we
apply Proposition~\ref{Sat-Order2-Stability-DDISRK} for RI6, then
we obtain $\mathcal{R}_{AS} = \{ (\hat{h}, k) \in \mathbb{C}^2 : |
(1+\hat{h} + \frac{1}{2} \hat{h}^2 + k^2)^2 - 3(1+\hat{h})^2 k^2 |
|1+\hat{h} + \frac{1}{2} \hat{h}^2 - \frac{1}{2} k^2|^4 < 1 \}$
and $\mathcal{R}_{MS} = \{ (\hat{h}, k) \in \mathbb{C}^2 : |1 +
\hat{h} + \frac{1}{2} \hat{h}^2 |^2 + |1 + \hat{h}|^2 |k|^2 +
\frac{1}{2} |k|^4 < 1\}$. The corresponding regions of stability
are given in Figure~\ref{Bild-RI6-DDIRDI4}.
Further, we can choose $c_1=0$, $c_2=\frac{1}{2}$ and e.g.\
$c_3=c_4=1$ which defines the scheme DDIRDI4. Then, the DDIRDI4
scheme (\ref{SRK-method-Ito-Wm-allg01}) is an advancement of the
stochastic theta method DDIRDI2 with $\theta=\frac{1}{2}$.
However, for DDIRDI4 we get $\mathcal{R}_{AS} = \{ (\hat{h}, k)
\in \mathbb{C}^2 : |(1+\frac{\hat{h}}{1-\frac{1}{2} \hat{h}} +
k^2)^2 - 3 (1+\frac{\hat{h}-\frac{1}{4} \hat{h}^2}{1-\frac{1}{2}
\hat{h}})^2 k^2| |1+ \frac{\hat{h}}{1-\frac{1}{2} \hat{h}} -
\frac{1}{2} k^2 |^4 < 1\}$ and $\mathcal{R}_{MS} = \{ (\hat{h}, k)
\in \mathbb{C}^2 : | 1 + \frac{\hat{h}}{1-\frac{1}{2} \hat{h}} |^2
+ |1 + \frac{\hat{h} - \frac{1}{4} \hat{h}^2}{1-\frac{1}{2}
\hat{h}}|^2 |k|^2 + \frac{1}{2} |k|^4 < 1\}$.
The regions of stability are given in
Figure~\ref{Bild-RI6-DDIRDI4}. Here, we can see that the good
stability properties of the order one scheme DDIRDI2 are not
carried over to the second order scheme DDIRDI4. Therefore, we are
looking for further second order DDISRK methods with some better
stability qualities.

\begin{figure}[tbp]
\begin{center}
\includegraphics[width=6.8cm]{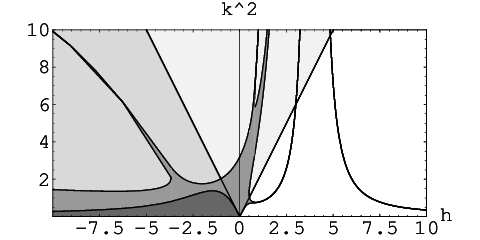}
\includegraphics[width=6.8cm]{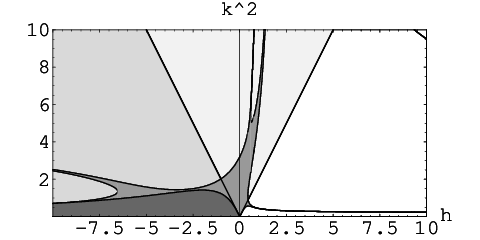}
\caption{Stability regions for DDIRDI5 with $c_1 = \frac{1}{4}$ on
the left and $c_1 = \frac{1}{2}$ on the right.}
\label{Bild-DDIRDI5-0.25-0.5}
\end{center}
\end{figure}

\begin{figure}[tbp]
\begin{center}
\includegraphics[width=6.8cm]{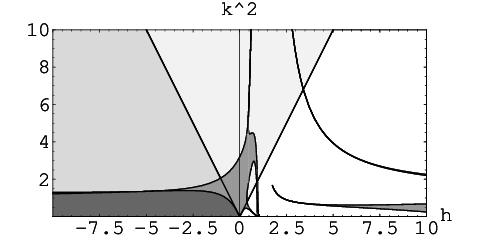}
\includegraphics[width=6.8cm]{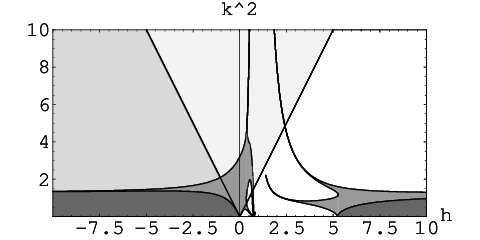}
\caption{Stability regions for DDIRDI5 with $c_1 = \frac{1}{2}+
\frac{\sqrt{3}}{6}$ on the left and $c_1 = 1$ on the right.}
\label{Bild-DDIRDI5-1.0-A}
\end{center}
\end{figure}

\begin{figure}[tbp]
\begin{center}
\includegraphics[width=6.8cm]{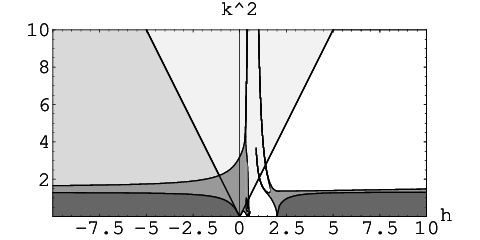}
\includegraphics[width=6.8cm]{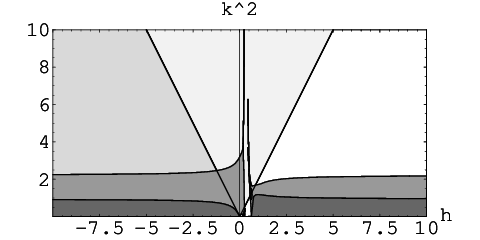}
\caption{Stability regions for DDIRDI5 with $c_1 = \frac{3}{2}$ on
the left and $c_1 = 3$ on the right.} \label{Bild-DDIRDI5-1.5-3.0}
\end{center}
\end{figure}

It is usual to consider singly diagonally implicit Runge--Kutta
methods for ODEs where all coefficients $A^{(0)}_{ii}$ are equal.
Therefore, we assume that $c_1=c_2$ for the schemes
(\ref{SRK-DDI-Coefficients}) in the following. Then, the domains
of stability are $\mathcal{R}_{AS} = \{ (\hat{h}, k) \in
\mathbb{C}^2 : |(1+ \frac{\hat{h} + (\frac{1}{2}-2 c_1)
\hat{h}^2}{(1-c_1 \hat{h})^2} + k^2)^2 - 3(1+ \frac{\hat{h} - c_1
\hat{h}^2}{(1-c_1 \hat{h})^2})^2 k^2 | |1+ \frac{\hat{h} +
(\frac{1}{2}-2 c_1) \hat{h}^2}{(1-c_1 \hat{h})^2} - \frac{1}{2}
k^2|^4 < 1\}$ and $\mathcal{R}_{MS} = \{ (\hat{h}, k) \in
\mathbb{C}^2 : |1+ \frac{\hat{h} + (\frac{1}{2}-2 c_1)
\hat{h}^2}{(1- c_1 \hat{h})^2} |^2 + |1+ \frac{\hat{h}- c_1
\hat{h}^2}{(1-c_1 \hat{h})^2}|^2 |k|^2 + \frac{1}{2} |k|^4 < 1 \}$
which depend on the coefficient $c_1$ of the scheme. In the
following, we consider various values for the parameter $c_1$ of
the DDISRK scheme (\ref{SRK-DDI-Coefficients}) and we analyze the
stability domain for $\lambda, \mu \in \mathbb{R}$. Therefore, we
choose $c_1 \in \{\frac{1}{4}, \frac{1}{2}, \frac{1}{2} +
\frac{1}{6} \sqrt{3}, 1, \frac{3}{2}, 3\}$. Especially, we
consider the case of $c_1 = c_2 = \frac{1}{2} + \frac{1}{6}
\sqrt{3}$ and $c_3=c_4=1$ which we denote as the scheme DDIRDI5.
Then, $A^{(0)}$ and $\alpha$ coincide with the coefficients of the
well known deterministic SDIRK scheme which is $A$-stable and
attains order $p_D=3$ for deterministic ODEs \cite{HW96}. The
corresponding stability regions are presented in
Figures~\ref{Bild-DDIRDI5-0.25-0.5}--\ref{Bild-DDIRDI5-1.5-3.0}.
\section{Numerical Experiments}
\label{Sec:Num-experiments}
We compare the efficiency of the proposed second order DDISRK
schemes DDIRDI4 and DDIRDI5 with the second order drift--implicit
SRK scheme DIPL1WM due to Platen~(\cite{KP99}, p.~501). Therefore,
we take the number of evaluations of the drift function $a$, of
the diffusion functions $b^j$, $j=1, \ldots, m$, and the number of
random numbers needed each step as a measure of the computational
complexity for each considered scheme.
\begin{figure}[tbp]
\begin{center}
\includegraphics[width=6.8cm]{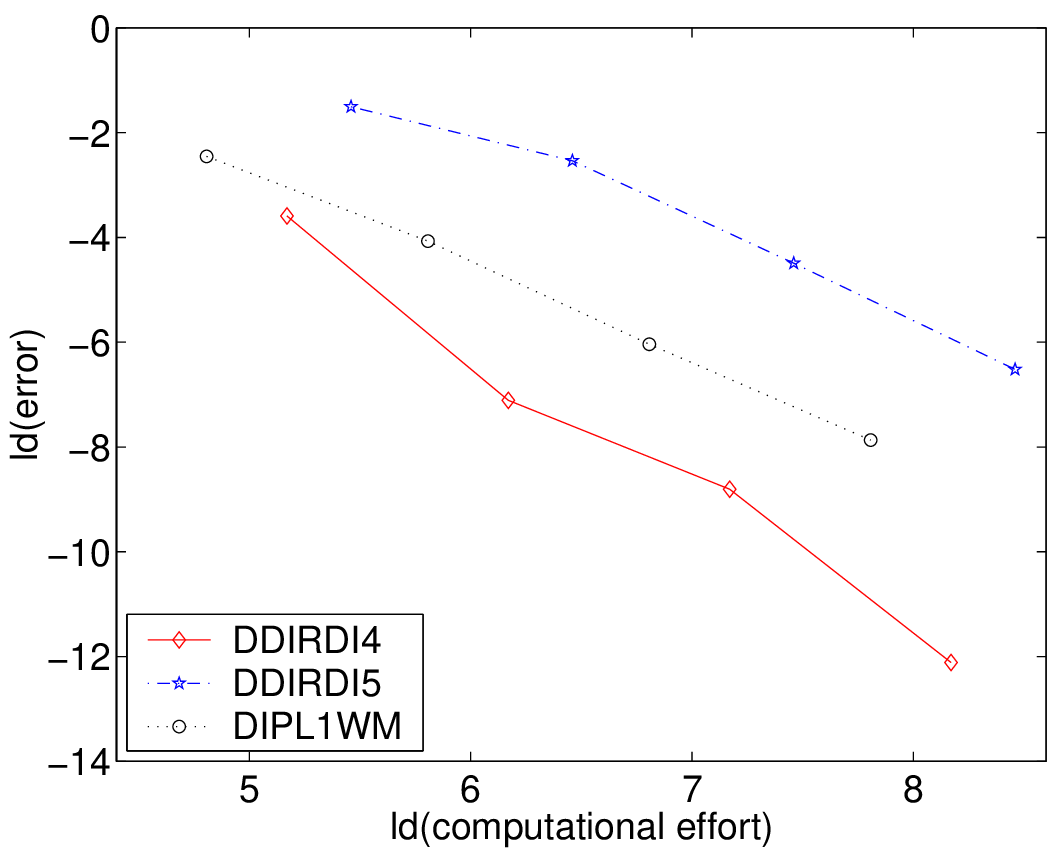}
\includegraphics[width=6.8cm]{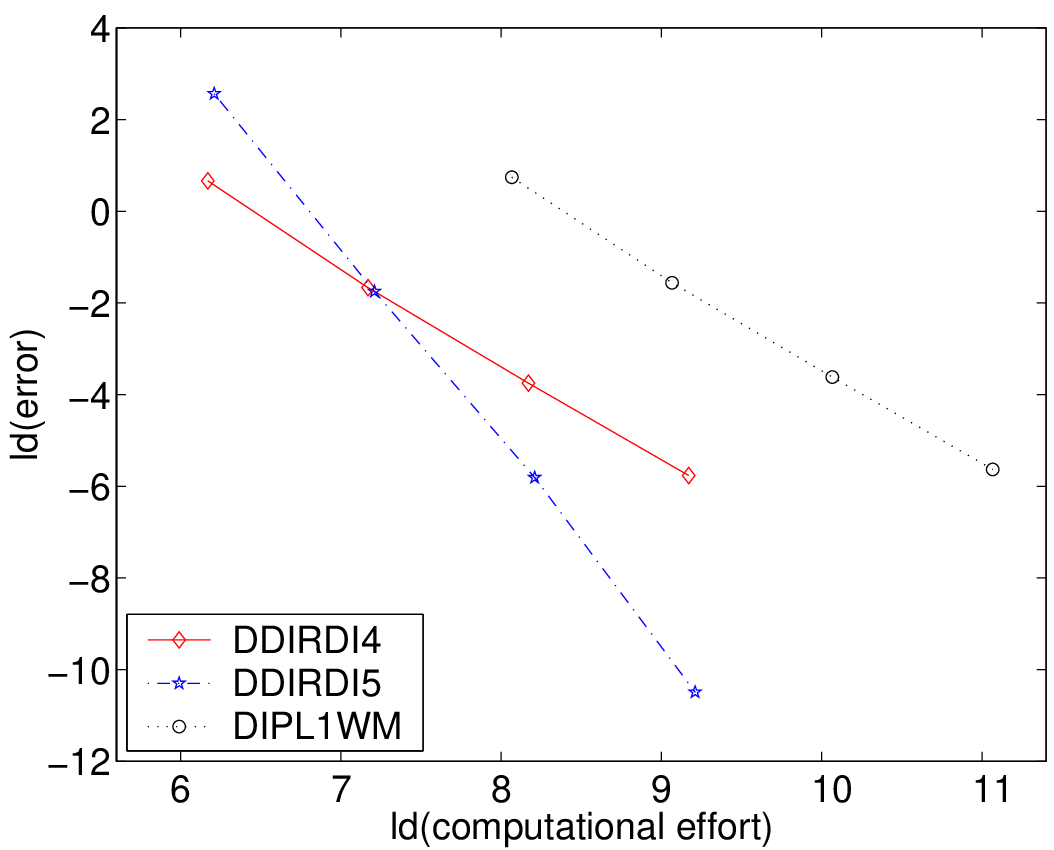}
\caption{Error vs computational effort with double logarithmic
scales for SDE (\ref{Simu:nonlinear-SDE2}) on the left and SDE
(\ref{Simu:nonlinear-SDE3}) on the right hand side.}
\label{Bild-effort}
\end{center}
\end{figure}
As a first example, we consider for $d=m=1$ the It\^{o} SDE
\begin{equation} \label{Simu:nonlinear-SDE2}
    {\mathrm{d}} X_t = \big( \tfrac{1}{2} X_t + \sqrt{X_t^2 + 1} \big) \,
    {\mathrm{d}}t + \sqrt{X_t^2 + 1} \, {\mathrm{d}}W_t \, , \qquad X_0=0,
\end{equation}
on the time interval $I=[0,2]$ with solution $X_t = \sinh (t +
W_t)$. Here, we choose $f(x)=p(\arsinh(x))$, where $p(z) = z^3 -
6z^2 + 8z$ is a polynomial. Then we calculate that $\E(f(X_t)) =
t^3 - 3t^2 + 2t$ which is approximated at time $t=2$ with step
sizes $2^{-1}, \ldots, 2^{-4}$ and $10^8$ simulated trajectories.
The results are presented on the left hand side of
Fig.~\ref{Bild-effort}.
As a second example, we consider a nonlinear SDE with a
$10$--dimensional driving Wiener process {\small{
\begin{equation} \label{Simu:nonlinear-SDE3}
    \begin{split}
    &{\mathrm{d}}X_t = X_t \, {\mathrm{d}}t
    + \frac{1}{10} \sqrt{X_t + \frac{1}{2}} \, {\mathrm{d}}W_t^1
    + \frac{1}{15} \sqrt{X_t + \frac{1}{4}} \, {\mathrm{d}}W_t^2
    + \frac{1}{20} \sqrt{X_t + \frac{1}{5}} \, {\mathrm{d}}W_t^3 \\
    &+ \frac{1}{25} \sqrt{X_t + \frac{1}{10}} \, {\mathrm{d}}W_t^4
    + \frac{1}{40} \sqrt{X_t + \frac{1}{20}} \, {\mathrm{d}}W_t^5
    + \frac{1}{25} \sqrt{X_t + \frac{1}{2}} \, {\mathrm{d}}W_t^6
    + \frac{1}{20} \sqrt{X_t + \frac{1}{4}} \, {\mathrm{d}}W_t^7 \\
    &+ \frac{1}{15} \sqrt{X_t + \frac{1}{5}} \, {\mathrm{d}}W_t^8
    + \frac{1}{20} \sqrt{X_t + \frac{1}{10}} \, {\mathrm{d}}W_t^9
    + \frac{1}{25} \sqrt{X_t + \frac{1}{20}} \, {\mathrm{d}}W_t^{10},
    \qquad X_0=1.
    \end{split}
\end{equation}
}} with non-commutative noise. Here, we approximate the second
moment of the solution $\E(X_t^2)= -\frac{68013}{14629060} +
(\frac{68013}{14629060}+1) \exp(\frac{731453}{360000} t)$ at time
$t=1$ by $10^8$ simulated trajectories with step sizes
$2^0,\dots,2^{-3}$. The results are presented on the right hand
side of Fig.~\ref{Bild-effort}. Here, the schemes DDIDRI4 and
DDIDRI5 perform impressively better than the drift--implicit
scheme DIPL1WM \cite{KP99}. This is a result of the reduced
complexity for the new class of efficient SRK schemes due to
R\"{o}{\ss}ler~\cite{Roe06c} which becomes significant especially for SDEs
with high-dimensional driving Wiener processes.

\begin{figure}[tbp]
\begin{center}
\includegraphics[width=6.8cm]{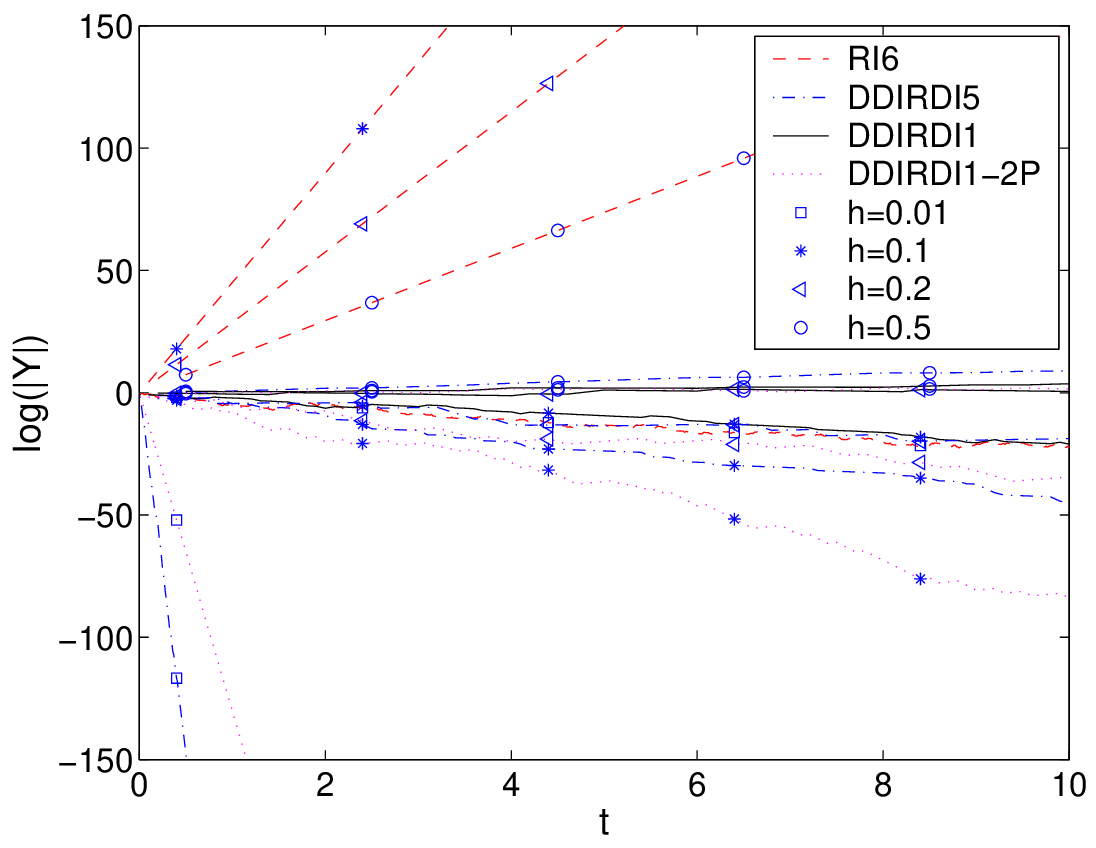}
\includegraphics[width=6.8cm]{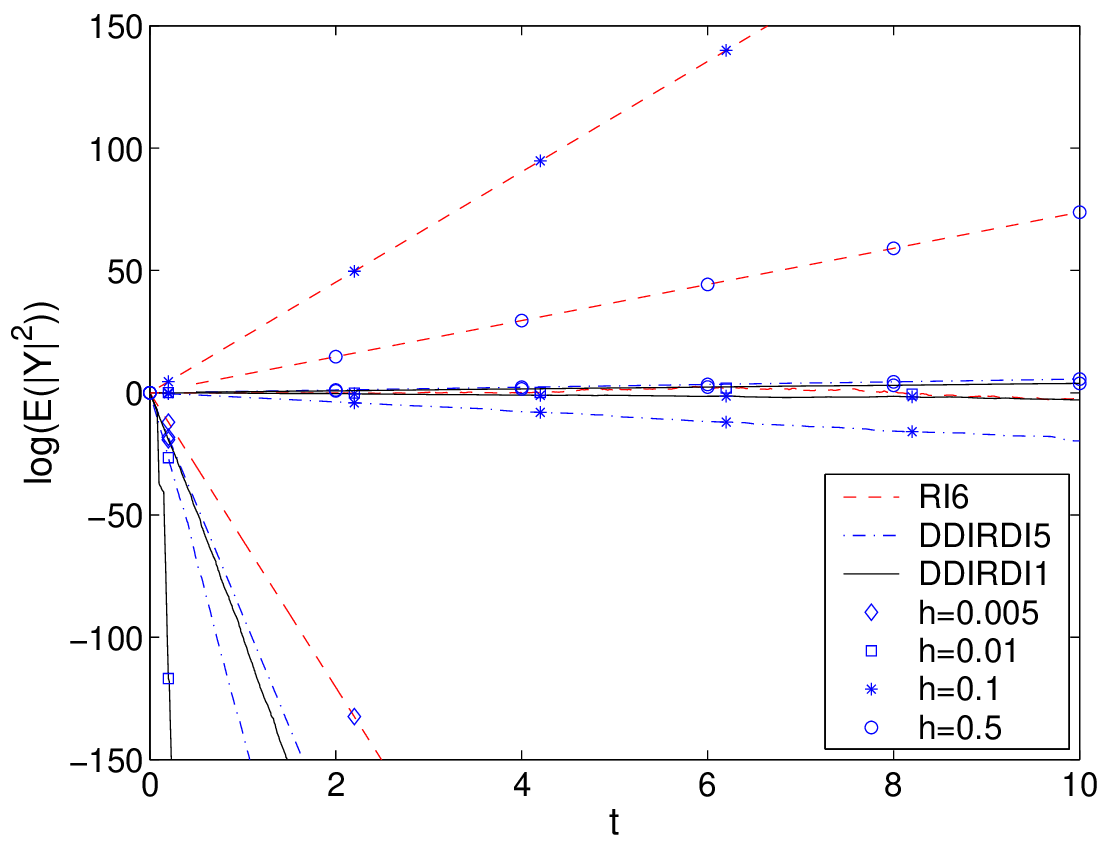}
\caption{Asymptotical and MS--stability analysis for RI6, DDIRDI1
and DDIRDI5.} \label{Bild-as-stability}
\end{center}
\end{figure}
Next, we verify the theoretical results for the domains of
stability of the proposed SRK methods by numerical experiments.
Therefore, we consider the test equation
(\ref{Lin-stoch-test-eqn}) with parameters $\lambda=-200$,
$\mu=\sqrt{5}$ and with initial value $X_0=1$ on the time interval
$[0,10]$. We apply the second order explicit SRK scheme RI6
\cite{Roe06c}, the order one DDISRK scheme DDIRDI1 with two point
as well as with three point distributed random variables and the
order two DDISRK scheme DDIRDI5. We denote by DDIRDI1-2P the
scheme DDIRDI1 if 2 point distributed random variables
$\tilde{I}_{(j),n}$ are used instead of $\hat{I}_{(j),n}$.

In order to analyse the numerical asymptotically stability, a
single approximation trajectory is simulated with each scheme
under consideration for the step sizes $h=0.01$, $h=0.1$, $h=0.2$
and $h=0.5$. Then, we obtain the following theoretical results due
to Proposition~\ref{Sat-1} and
Proposition~\ref{Sat-Order2-Stability-DDISRK}: the scheme RI6 is
asymptotical stable for the step size $h=0.01$ and it is unstable
for $h=0.1$, $h=0.2$ and $h=0.5$. DDIRDI1 and DDIRDI1-2P are
stable for $h=0.01$ and $h=0.1$. In the case of $h=0.2$ only
DDIRDI1-2P is stable while DDIRDI1 is unstable. Further, DDIRDI1
and DDIRDI1-2P are unstable for $h=0.5$. Finally, the scheme
DDIRDI5 is asymptotical stable for $h=0.01$, $h=0.1$ and even for
$h=0.2$, however it is unstable for $h=0.5$. The numerical results
for a single trajectory $|Y_n|$ are plotted with logarithmic scale
to the base 10 versus the time on the left hand side of
Fig.~\ref{Bild-as-stability}. We remark that the results for
DDIRDI1 with step size $h=0.01$ tend to zero after two steps and
are thus not visible in Fig.~\ref{Bild-as-stability}.

For the analysis of the numerical MS--stability, the value
$\E(|X_t|^2)$ is approximated by Monte Carlo simulation with
$10^{4}$ independent trajectories for the step sizes $h=0.005$,
$h=0.01$, $h=0.1$ and $h=0.5$. Proposition~\ref{Sat-1} and
Proposition~\ref{Sat-Order2-Stability-DDISRK} give the following
results: RI6 is MS--stable for $h=0.005$ and MS--unstable for all
other considered step sizes. DDIRDI1 and DDIRDI1-2P are MS--stable
for $h=0.005$ and $h=0.01$, however MS--unstable for $h=0.1$ and
$h=0.5$. Further, DDIRDI5 is MS--stable for step sizes $h=0.005$,
$h=0.01$ and even for $h=0.1$ and MS--unstable for $h=0.5$. The
corresponding numerical results of $\E(|Y_n|^2)$ are presented
with logarithmic scale to the base 10 versus the time on the right
hand side of Fig.~\ref{Bild-as-stability}. Again, the numerical
results exactly confirm our theoretical findings for the domains
of stability.

\end{document}